\documentclass{amsart}
\usepackage{amssymb,amsmath,latexsym}
\usepackage{amsthm}
\usepackage{fontenc}
\usepackage{amssymb}

\numberwithin{equation}{section}

\newtheorem{theorem}{Theorem}[section]

\begin{document}
\author{Alexander E. Patkowski}
\title{A general size-biased distribution}

\maketitle

\begin{abstract} We generalize a size-biased distribution related to the Riemann xi function using the work of Ferrar. Some analysis and properties of this more general distribution are offered as well.\end{abstract}



\section{Introduction} 
 The Riemann xi function is defined to be [4, 11] $\xi(s):=\frac{1}{2}s(s-1)\pi^{-\frac{s}{2}}\Gamma(\frac{s}{2})\zeta(s),$ where $\zeta(s)$ is the Riemann zeta function. It is a well-known property that $\xi(s)=\xi(1-s).$ The integral representation that provides the probabilistic density function contained in [1] is found in [4, pg.207--208],
\begin{equation}\int_{0}^{\infty}x^{s-1}\Theta(x)dx=\xi(s),\end{equation}
where 
\begin{equation}\Theta(x):=2x^2\sum_{n\ge1}(2\pi^2 n^4x^2-3\pi n^2)e^{-\pi n^2 x^2},\end{equation} 
valid for $s\in \mathbb{C}.$
\par The function $\Theta(x)$ also enjoys the nice functional property that $x\Theta(x)=\Theta(x^{-1}).$ It is this property that gives rise to its sized-biased distribution property, as can be seen in [1, 3]. Namely, we have that $\mathop{\mathbb{E}(Xf(X))}=\mathop{\mathbb{E}(f(\frac{1}{X}))},$ for a measurable function $f,$ and random variable $X$ with density $\frac{1}{x}\Theta(x).$ See [9] for more information on this distribution and applications. In particular, it can be found in [1, 3] that $2\xi(0)=1,$ and hence also $2\xi(1)=1,$ which suggests one property [2, pg.34, F2] of a probability distribution (on $(0,\infty)$) given (1.1). In addition, it has been discovered that, for a random variable $X$ with density function $x^{-1}\Theta(x)$ [1, 3],
$$P(X\le x)=-x^{-2}\Theta'(x^{-1}),$$
$$\mathop{\mathbb{E}(X^s)}=2\xi(s),$$
and \begin{equation}P(X\le x)=4\pi x^{-3}\sum_{n\ge1}n^2e^{-\pi n^2/x^2}.\end{equation} Here we denote $\mathop{\mathbb{E}}$ and $\mathop{\mathbb{V}}$ to be the expectation and variance, respectively. For more on this distribution see also [10].
\par In [6] Ferrar made use of the underlying functional relationship $\xi(s)=\xi(1-s)$ to establish that the function 
\begin{equation} w_k(x):=\frac{1}{2\pi i}\int_{(c)}(\Gamma(\frac{s}{2})\zeta(s)\pi^{-s/2})^kx^{-s}ds, \end{equation}
for integers $k\ge1,$ and real $c>1,$ satisfies $w_k(x)-R(x,k)=-\frac{1}{x}R(\frac{1}{x},k)+\frac{1}{x}w_k(\frac{1}{x}).$ Here we used Ferrar's notation to denote $R(x,k)$ as the residue at the pole $s=0,$ corresponding to $k.$ For $k=1$ we find a connection the the Riemann xi function and consequently the probability distribution we have noted. \par The purpose of this paper is to generalize the distribution of the Riemann xi function using the work of Ferrar and offer some of its properties. In doing so we provide a new general size-biased distribution. \par Define the differential operator by:
$$D^{1}_x=\frac{\partial }{\partial x}(x^2\frac{\partial }{\partial x}),$$
$$D^{2}_x=\frac{\partial }{\partial x}(x^2\frac{\partial^2 }{\partial x^2}x^2\frac{\partial }{\partial x}),$$
$$D^{k}_x=\frac{\partial }{\partial x}x^2\underbrace{\frac{\partial^2 }{\partial x^2}x^2\cdots\frac{\partial^2 }{\partial x^2}x^2}_{k-1}\frac{\partial }{\partial x}.$$
Throughout we will define $$v_k(x):=D^{k}_x (w_k(x)).$$ Note that by [1, pg.439] or [3], $v_1(x)=\Theta(x)>0$ when $x>0.$ It is a simple observation using Parseval's theorem for Mellin transforms and the functional equation for $\xi(s)$ to see that $v_2(x)>0$ when $x>0.$ Specifically, the resulting integral $$v_2(x)=\int_{0}^{\infty}\Theta(y)\Theta(yx)dy,$$ must also be positive for $x>0$ since $v_1(x)>0.$ Repeating this similarly, we see that 
$$v_3(x)=\int_{0}^{\infty}\Theta(zx)\int_{0}^{\infty}\Theta(y)\Theta(yz)dydz>0,$$ for $x>0.$ Iterating this idea repeatedly proves $v_k(x)>0$ for $k\ge1.$
\begin{theorem} Let $X_k$ be a random variable on $(0,\infty),$ with density function $\frac{1}{x}v_k(x).$ We have, $v_k(x)=\frac{1}{x}v_k(\frac{1}{x}).$ Then, consequently, $1/X_k$ is equivalent to the size-biased distribution from $X_k$ and hence $$\mathop{\mathbb{E}(X_kf(X_k))}=\mathop{\mathbb{E}(f(\frac{1}{X_k}))},$$
and
$$\mathop{\mathbb{E}(X_k^s)}=\left(2\xi(s)\right)^{k}.$$ Additionally, we have $\mathop{\mathbb{E}(X_k)}=1,$ and $\mathop{\mathbb{V}(X_k)}=(2\zeta(2))^k-1.$
\end{theorem}
Ferrar [6] had written the function explicitly in terms of $d_k(n),$ the number of ways of representing $n$ as a product of $k$ factors. However, there is no simple known formula for $$\frac{1}{2\pi i}\int_{(c)}(\Gamma(\frac{s}{2})\pi^{-s/2})^kx^{-s}ds,$$ when $k\ge3.$ Although, it is possible to work with Parseval's theorem to create nested integrals as an alternative expression. For $k=2,$ the integral may be shown to be connected to $K_m(x),$ the modified Bessel function of the second kind for $m=0.$ Since $\frac{\partial }{\partial x}K_0(ax)=-aK_1(ax),$ and $\frac{\partial }{\partial x}K_1(ax)=-a(K_0(ax)+\frac{1}{x}K_1(ax)),$ we have $$v_2(x)=\sum_{n\ge1}d_2(n)(2\pi n x)^2\left(((2\pi n x)^2+9)K_0(2\pi n x)-6(2\pi n x)K_1(2\pi n x)\right).$$ Convergence of this series is absolute since $K_m(x)$ decays exponentially as $x\rightarrow\infty$ [8, pg.250]. We are able to prove that $v_2(x)$ is non-negative for $x\ge 1.128$ by the first term of the series alone. First, recall [7, pg.385, (A.9), (A.11)] that $K_2(x)=\frac{2}{x}K_1(x)+K_0(x),$ and $K_{m}(x)\ge K_{m-1}(x)$ for $x>0,$ $v\ge\frac{1}{2}.$ Hence $K_0(x)\ge (1-\frac{2}{x})K_1(x).$ Therefore,
$$\begin{aligned} &((2\pi n x)^2+9)K_0(2\pi n x)-6(2\pi n x)K_1(2\pi n x)\\
& \ge \left((2\pi n x)^2+9\right)\left(1-\frac{1}{\pi n x}\right)K_1(2\pi n x)-6(2\pi n x)K_1(2\pi n x)\\
&\left((2\pi n x)^2+9-4n\pi x -\frac{9}{\pi n x}-6(2\pi n x)\right)K_1(2\pi n x). \end{aligned}$$
Since $K_m(x)>0,$ for $x>0,$ and $m$ any real number [7, pg.384], it suffices to check the polynomial in the last line. If $n=1,$ this factor is non-negative if $x<0$ or $x\ge1.128.$ Inspection of the terms $n>1$ show that each are non-negative for $x<0,$ $x\ge 1.128,$ and certain ranges contained in $0<x< 1.128.$ That is, the region for which the terms of the sum are negative become smaller as $n$ becomes larger.
\begin{theorem} The distribution $X_{k-1}$ majorises $X_k$ in distribution, for each $k>1.$ \end{theorem}
We next give an interesting transformation property giving a connection between the $k=2$ and $k=1$ distributions through series identities. 
\begin{theorem} For $x>0$ we have,
$$x\sum_{n\ge1}d_2(n)nK_1(2xn\sqrt{\pi})=\frac{1}{4\sqrt{\pi}}\int_{0}^{\infty}P(X_1\le y)\left(\sum_{n\ge1}e^{-(xny)^2}\right)dy. $$
\end{theorem}
Recall the notation [2, pg.726] for convergence in probability  $\overset{p}{\to}.$ Define $S_n:=\sum_{i=1}^{n}X_{k_i}.$ Now we are able apply Khintchin's law of large numbers to Theorem 1.1.
\begin{theorem} Let $\{X_{k_i}\}$ be a sequence (in $i$) of independent random variables all distributed as defined in Theorem 1.1. Let $S_n=\sum_{i\ge1}^{n}X_{k_i}.$ If $S_n/n$ converges in probability, then, $$\frac{S_n}{n}\overset{p}{\to}1,$$ as $n\rightarrow\infty.$ Equivalently, 
$$\mathop{\mathbb{E}\left|\frac{S_n}{n}-1\right|}\rightarrow0,$$
as $n\rightarrow\infty.$
\end{theorem}

\section{The proofs of Main Results}
To prove our theorem we will use standard properties of Mellin transforms, and refer the reader to [9] for an introduction.
\begin{proof}[Proof of Theorem 1.1] Taking Mellin transform of $v_k(x),$ we have
$$\int_{0}^{\infty} x^{s-1}v_k(x)dx=(s(s-1)\Gamma(\frac{s}{2})\zeta(s)\pi^{-s/2})^k,$$
for $\Re(s)>1.$ It follows from analytic continuation, similar to [4, 10], that we also have $s\in\mathbb{C}.$ One way to see this is to observe that the residue at the poles $s=0$ and $s=1$ of $(s(s-1)\Gamma(\frac{s}{2})\zeta(s)\pi^{-s/2})^k$ is $0.$ Hence we may take the integral $v_k(x),$ compute the residue at these poles moving the line of integration to $-1<\Re(s)<0,$ then replace $s$ by $1-s$ returning the integral to the region $\Re(s)=c>1.$ It follows that replacing $x$ by $1/x,$ and then applying the functional equation of Ferrar, we find that we have the formula $\mathop{\mathbb{E}(X_k^s)}=2^k\xi^{k}(s).$ The sized-biased property follows from the simple change of variable $x$ to $\frac{1}{x}$ and then applying the functional equation $v_k(x)=\frac{1}{x}v_k(\frac{1}{x}).$ The variance value is computed with $\mathop{\mathbb{V}(X_k)}=\mathop{\mathbb{E}(X_k^2)}-\mathop{\mathbb{E}(X_k)^2}.$ \end{proof}
\begin{proof}[Proof of Theorem 1.2] Using the definition given in [2, pg. 302], we need to prove that
\begin{equation} P(X_{k}\ge x)\le P(X_{k-1}\ge x). \end{equation}
Therefore, we need to only show that 
$$\int_{x}^{\infty}y^{-1}v_k(y)dy\le\int_{x}^{\infty}y^{-1}v_{k-1}(y)dy.$$ 
Now directly integrating shows that for $x>0,$ $c>0,$
$$\begin{aligned} &\int_{x}^{\infty}y^{-1}v_k(y)dy= \int_{x}^{\infty}y^{-1}\left(\frac{1}{2\pi i}\int_{(c)}(s(s-1)\Gamma(\frac{s}{2})\zeta(s)\pi^{-s/2})^ky^{-s}ds\right)dy \\
&=-\frac{1}{2\pi i}\int_{(c)}(s(s-1)\Gamma(\frac{s}{2})\zeta(s)\pi^{-s/2})^k\frac{x^{-s}}{s}ds, \end{aligned}$$
by an application of Fubini's theorem. Moreover, using Stirling's formula and arguments found in [9, pg.121] for (1.4), we see that $v_{k-1}(x)\gg v_k(x),$ for $x>0.$ Since $v_k(x)>0,$ for $x>0$ as we showed previously, it follows that $v_{k-1}(x)\geq v_k(x),$ for $x>0.$
Parseval's theorem for Mellin transforms [9] may be applied with (1.1)--(1.2) to show that 
$$\int_{0}^{\infty}y^{-1}v_k(y)dy=\int_{0}^{\infty}y^{-1}v_{k-1}(y)2y^2\sum_{n\ge1}(2\pi^2 n^4y^2-3\pi n^2)e^{-\pi n^2 y^2}dy,$$ since the common region of holomorphy is $\mathbb{C}.$ The latter integral ensures that the integrand of the first integral decays to $0$ as $y\rightarrow\infty$ faster than just $y^{-1}v_{k-1}(y).$ This covers the $x\rightarrow0^{+}$ case. Hence (2.1) is established. \end{proof}
\begin{proof}[Proof of Theorem 1.3] We start with an integral that may be found in [5, pg.313, $h=1$] and making a change of variable $y=1/y$,
\begin{equation}2K_{s}(2\alpha\beta)=(\frac{\alpha}{\beta})^{s}\int_{0}^{\infty}e^{-\beta^2y^2-\frac{\alpha^2}{y^2}}\frac{dy}{y^{2s+1}},\end{equation}
valid for $\Re(\alpha)>0,$ $\Re(\beta)>0.$ We select $s=1$ and put $\alpha=\sqrt{\pi}n$ in (2.2) and sum over $n$ to obtain for $x>0$
\begin{equation}x\sum_{n\ge1}nK_1(2xn\sqrt{\pi})=\sqrt{\pi}\int_{0}^{\infty}e^{-(xy)^2}\left(\frac{1}{y^3}\sum_{n\ge1}n^2e^{-\pi n^2/y^2}\right)dy.\end{equation} The exchange of integral and series is justified by absolute convergence. Replacing $x$ by $mx$ in (2.3) and again summing over $m$ gives the result when comparing with (1.3).\end{proof}
\begin{proof}[Proof of Theorem 1.4] This theorem is a direct application of [2, pg.185, Theorem 8.1.1] with the property $\mathbb{E}(X_k)=1,$ and includes the equivalent statement [2, pg.186, Corollary 8.1.1].

\end{proof}

1390 Bumps River Rd. \\*
Centerville, MA
02632 \\*
USA \\*
E-mail: alexpatk@hotmail.com, alexepatkowski@gmail.com

\end{document}